\documentclass[12pt]{article}

\date{} 

\begin{document}

\centerline{} 

\centerline{} 

\centerline {\large{\bf Discrete duality for Tense Symmetric Heyting Algebras}} 

\centerline{} 


\centerline{} 

\centerline{\bf {Aldo V. Figallo}} 

\centerline{} 

\centerline{Universidad Nacional de San Juan. Instituto de Ciencias B\'asicas} 

\centerline{Avda. I. de la Roza 230 (O), 5400 San Juan, Argentina} 

\centerline{and} 

\centerline{Universidad Nacional del Sur. Departamento de Matem\'atica} 

\centerline{Avda. Alem 1253, 8000 Bah\'ia Blanca, Argentina} 

\centerline{}

\centerline{\bf {Gustavo Pelaitay}} 

\centerline{} 

\centerline{Universidad Nacional de San Juan. Instituto de Ciencias B\'asicas} 

\centerline{Avda. I. de la Roza 230 (O), 5400 San Juan, Argentina} 

\centerline{and} 

\centerline{Universidad Nacional del Sur. Departamento de Matem\'atica} 

\centerline{Avda. Alem 1253, 8000 Bah\'ia Blanca, Argentina} 

\centerline{gpelaitay@gmail.com}

\centerline{} 

\centerline{\bf {Claudia Sanza}} 

\centerline{} 

\centerline{Universidad Nacional del Sur. Departamento de Matem\'atica} 

\centerline{Avda. Alem 1253, 8000 Bah\'ia Blanca, Argentina} 

\centerline{}

\newtheorem{theorem}{\quad Theorem}[section] 

\newtheorem{proof}{\quad Proof} 

\newtheorem{definition}{\quad Definition}

\newtheorem{lemma}{\quad Lemma} 

\newtheorem{remark}{\quad Remark} 

\newtheorem{Example}{\quad Example}

\newtheorem{proposition}{\quad Proposition}

\begin{abstract}In this article, we continue the study of tense symmetric Heyting algebras (or $TSH$-algebras). These algebras constitute a generalization of tense algebras. In particular, we describe a discrete duality for $TSH$-algebras bearing in mind the results indicated by E. Or\l owska and I. Rewitzky in [E. Or\l owska and I. Rewitzky,
           {\it Discrete Dualities for Heyting Algebras with Operators}, 
            Fund. Inform. {\bf 81} (2007), no.1--3, 275--295.] for Heyting algebras.
In addition, we introduce a propositional calculus and prove this calculus has $TSH$-algebras as algebraic counterpart. Finally, the duality mentioned above allowed us to show the completeness theorem for this calculus.
\end{abstract} 



\section{\large Introduction and preliminaries}

A discrete duality (see \cite{EO3, EO, EO1}) is a duality between classes of algebras and classes of relational
systems (frames):

Let {\bf Alg} be a class of algebras and let {\bf Frm} be a class of frames.

Establishing a discrete duality between these two classes requires the following steps:

\begin{enumerate}
\item[(i)] With every algebra $W$ from {\bf Alg} associate a canonical frame ${\mathcal X}(W)$ of the algebra
and show that it belongs to {\bf Frm}.
\item[(ii)] With every frame $X$ from {\bf Frm} associate a complex algebra ${\mathcal C}(X)$, and show
that it belongs to {\bf Alg}.
\item[(iii)] Prove two Representation Theorems:
\begin{enumerate} 
\item For each $W\in {\bf Alg}$ there is an embedding $h:W\to {\mathcal C}({\mathcal X}(W))$.
\item For each frame $X\in {\bf Frm}$ there is an embedding $k:X\to {\mathcal X}({\mathcal C}(X))$.
\end{enumerate}
\item[]
\end{enumerate} 
An important application of discrete duality is that it provides a Kripke semantics (resp. an algebraic semantics) once an algebraic semantics (resp. a Kripke semantics) for a formal language is given (see \cite{EO1}).

\vspace{1mm}

Let $T$ be a binary relation on a set $X$ and let $A$ be a subset of $X$. In what follows we will denote by $[T]A$ the set $\{x\in X:\, {\rm for\,\, all}\,\, y, x\, T\, y\,\, {\rm implies}\,\, y\in A\}$.

\vspace{1mm}

In \cite{EO}, E. Orlowska and I. Rewitzky introduced the notion of Heyting frame (or $H$-frame, for short) as a pair $(X, \leq)$ where $X$ is a non-empty set and $\leq$ is a quasi-order on $X$. These authors proved that if  $\langle W, \vee, \wedge, \to, 0, 1\rangle$ is a Heyting algebra, then its canonical frame is $({\mathcal X}(W),\leq^{c})$, where  ${\mathcal X}(W)$ is the set of all prime filters of $W$ and  $\leq^{c}$ is $\subseteq$. It is easy to see that this canonical frame
is an $H$-frame. On the other hand, given an $H$-frame $(X, \leq)$, they show that its complex algebra is $\langle{\mathcal C}(X), \vee^{c},\wedge^{c}, \to^{c}, 0^{c}, 1^{c}\rangle$, where ${\mathcal C}(X)=\{A\subseteq X: [\leq]A=A\}$, $0^{c}=\emptyset$, $1^{c}=X$, $A\vee^{c}B=A\cup B$, $A\wedge^{c}B=A\cap B$ and  $A\to^{c}B=[\leq]((X\setminus A)\cup B)$ for all $A,B\in{\mathcal C}(X)$.

\vspace{1mm}

These results allowed them to obtain a discrete duality for Heyting algebras by defining the embeddings as follows:
\begin{itemize}
\item[(E1)] $h:W \rightarrow {\mathcal C}({\mathcal X}(W))$, $h(a)=\{F\in {\mathcal X}(W): a\in F\}$, 
\item[(E2)] $k: X\rightarrow {\mathcal X}({\mathcal C}(X))$, $k(x)=\{A\in {\mathcal C}(X): x\in A\}$. 
\end{itemize}

\vspace{1mm}

On the other hand, in 1942 Gr. C. Moisil (\cite{Mo}) introduced the modal symmetric propositional calculus as an extension of the positive calculus of Hilbert-Bernays obtained by adding a new negation connective, $\sim$, the axiom schemata
$$\alpha\to\sim\sim \alpha, \hspace{1cm}\sim\sim \alpha\to \alpha$$

\noindent and the contraposition rule

\begin{center}
if $\alpha\to\beta $ then $\sim \beta\to \sim \alpha$.
\end{center}

This propositional calculus has symmetric Heyting algebras as the algebraic counterpart. These algebras were investigated by A. Monteiro (\cite{M}) and also by L. Iturrioz (\cite{I2}) and H.P. Sankappanavar (\cite{S}).

\vspace{1mm}

Recall that an algebra $\langle W,\vee,\wedge,\to,\sim,0,1\rangle$ is a symmetric Heyting algebra (see \cite{M}) if $\langle W,\vee,\wedge,\sim,0,1\rangle$ is a De Morgan algebra and $\langle W,\vee,\wedge,\to,0,1\rangle$ is a Heyting algebra.

\vspace{1mm}

In \cite{FPS}, we introduce tense symmetric Heyting algebras (or $TSH$-algebras, for short) as algebras $\langle W,\vee,\wedge,\to,\sim,G,H,0,1\rangle$ where the reduct $\langle W,\vee,\wedge,\to,\sim,0,1\rangle$ is a symmetric  Heyting algebra  and $G,\,H$ are unary operators on $W$ verifying  these conditions: 
 
\begin{itemize}
\item[\rm (T1)] $G(1)=1$, $H(1)=1$,
\item[\rm (T2)] $G(x\wedge y)=G(x)\wedge G(y)$, $H(x\wedge y)=H(x)\wedge H(y)$,
\item[\rm (T3)] $x\leq G (\sim H (\sim x))$,  $x\leq H(\sim G (\sim x)) $.
\end{itemize}

In what follows, we will denote these algebras by $(W,G,H)$ or simply by $W$.

\begin{remark}\label{rem.3.2.}

\noindent (i)\, From \textnormal{(T2)} it follows that $G$ and $H$ are increasing.

\noindent (ii)\, If $\langle A,\vee,\wedge,\to,\sim,G,H,0,1\rangle$ is a $TSH$--algebra in which every element of $A$ is boolean, then $\langle A,\vee,\wedge,\sim,G,H,0,1\rangle$ is a tense algebra (\cite{K,PU}).
\end{remark}

\section{\large A discrete duality for $TSH$-algebras}

In this section, we describe a discrete duality for $TSH$-algebras taking into account the one indicated above for Heyting algebras. To this end, we introduce the following

\begin{definition} A $TSH$-frame is a structure $(X,\leq, g, R, Q)$ where $(X,\leq)$ is a $H$-frame, $g:X\rightarrow X$ is a function, $R,Q$ are binary relations on $X$ and the following conditions are satisfied:
\begin{itemize}
\item [\rm (K1)] if $x\leq y$ then $g(y)\leq g(x)$ for $x, y \in X$,
\item [\rm (K2)] $g(g(x))=x$ for $x\in X$,
\item [\rm (K3)] $(\leq\circ R\circ\leq)\subseteq R$,
\item [\rm (K4)] $(\leq\circ Q\circ\leq)\subseteq Q$,
\item [\rm (K5)] $x\,R\,g(y)$ if and only if $y\,Q\,g(x)$ for $x, y \in X$. 
\end{itemize}
\end{definition}

In what follows, $TSH$-frames will be denoted simply by $X$ when no confusion may arise.

\begin{definition}
A canonical frame of a $TSH$-algebra $(W,G,H)$ is a structure $({\mathcal X}(W),\leq^{c},g^{c},R^{c},Q^{c})$, where $({\mathcal X}(W),\leq^{c})$ is the canonical frame associated with $\langle W,\vee,\wedge,\to,0,1\rangle$ and the following conditions are verified for $P, F\in{\mathcal X}(W)$:
\begin{itemize}
\item[\rm (F1)] $g^{c}(P)=\{a\in W:\, \sim a\notin P\}$,
\item[\rm (F2)] $PR^{c}F$ if and only if $G^{-1}(P)\subseteq F$, 
\item[\rm (F3)] $PQ^{c}F$ if and only if $H^{-1}(P)\subseteq F$.
\end{itemize}
\end{definition}

\vspace{1mm}

\begin{lemma}
The canonical frame of a $TSH$-algebra is a $TSH$-frame.
\end{lemma}
{\bf Proof} Taking into account the results established in \cite[Lemma 11.1]{EO3}, we only have to prove  (K3), (K4) and (K5).

\noindent (K3): Let $(P,F)\in\, \leq ^{c}\circ R^{c}\circ \leq^{c}$. Then there exist $T,S\in{\mathcal X}(W)$ such that  $P\subseteq T$, $TR^{c}S$ and  $S\subseteq F.$ From the last two assertions we have that $G^{-1}(T)\subseteq F$. Therefore,  since $P\subseteq T$ we infer  that $P\,R^{c}\, F$. 

\noindent (K4): It is proved in a similar way to (K3).

\noindent (K5): Let $F\, R^{c}g^{c}(P)$ and $a\in H^{-1}(P)$. Suppose that $\sim a\in F$. On the other hand, from (T3) we have that  $\sim a\leq G(\sim H(a))$ and so, we get that $G(\sim H(a))\in F$. From this last assertion and the fact that $G^{-1}(F)\subseteq g^{c}(P)$, we obtain $\sim H(a)\in g^{c}(P)$.  Hence, $H(a)\notin P$ which is a contradiction. Therefore, $a\in g^{c}(F)$ from which we conclude that $PQ^{c}g^C(F)$. The converse is proved similarly.

\begin{definition}The complex algebra of a $TSH$-frame $(X,\leq,g,R,Q)$ is\\ $\langle {\mathcal C}(X),\vee^{c},\wedge^{c},\to^{c},\sim^{c},G^{c},H^{c},0^{c},1^{c}\rangle$, where $\langle {\mathcal C}(X),\vee^{c},\wedge^{c},\to^{c},0^{c},1^{c}\rangle$ is the complex algebra of the $H$-frame $(X,\leq)$, $\sim^{c}A=X\setminus g(A)$, $G^{c}(A)=[R]A$ and $H^{c}(A)=[Q]A$, for all $A\in{\mathcal C}(X)$.
\end{definition}

\vspace{1mm}

\begin{lemma} The complex algebra of a $TSH$-frame is a $TSH$-algebra.
\end{lemma}
{\bf Proof} From \cite{EO3, EO}, ${\mathcal C}(X)$ is closed under the lattice operations, $\sim^{c}$ and  $\to^{c}$. Now, we show that it is also  closed under $G^{c}$ i.e., $G^{c}A=[\leq]G^{c}A$. From the reflexivity of $\leq$, we have that $[\leq]G^{c}A \subseteq G^{c}A$. Assume that  $x\in G^{c}A$. Let $y\in X$ be such that $x\leq y$ and take any $z\in X$ verifying $yRz$. Hence, from the reflexivity of $\leq$ and (K3) we infer that $xRz$. So, $z\in A$ and therefore, $x\in [\leq]G^{c}A$. Thus, $G^{c}A\subseteq[\leq]G^{c}A$. Similarly, it is proved that $H^{c}A=[\leq]H^{c}A$. On the other hand, clearly (T1) and (T2) are verified. Therefore, it only remains to prove (T3). Let $x\in A$ and suppose that $x\notin G^c(\sim^{c}H^c(\sim^cA))$. Then there is $y$ such that $xRy$ and $y\notin \sim^{c}H^c(\sim^{c}A)$. From this last statement, $y\in g(H^c(\sim^{c}A))$ and so,  $y=g(z)$ for some $z\in H^c(\sim^{c}A)$. Hence, $xRg(z)$ and from (K5) we get that  $z Q g(x)$. This assertion and the fact that $z\in H^c(\sim^{c}A)$ enable us to infer that $g(x)\notin g(A)$, which is a contradiction. So, $A\subseteq G^c(\sim^{c}H^c(\sim^cA))$. Analogously, it is proved that $A\subseteq H^c(\sim^{c}G^c(\sim^cA))$.

\begin{remark}
Let $(W,G,H)$ be a $TSH$-algebra. It is worth noting that if $F$ is a filter of $W$ then $G^{-1}(F)$ and $H^{-1}(F)$ are also filters of $W$.
\end{remark}

\vspace{1mm}

\begin{theorem}\label{t21} 
Each $TSH$-algebra $W$ is embeddable into ${\mathcal C}({\mathcal X}(W))$.
\end{theorem}

{\bf Proof} Let us consider the function $h : W \to{\mathcal C}({\mathcal X}(W))$ defined by $h(a)=\{P\in{\mathcal X}(W): a\in P\}$,
for all $a\in W$ (see \cite{EO3, EO}).  Let $F\in h(G(a))$;  then $G(a)\in F$. Suppose that $P\in {\mathcal X}(W)$ verifies that $FR^cP$. Then from (F2), $G^{-1}(F)\subseteq P$ and so,  $a\in P$. Therefore, $F\in G^{c}(h(a))$ from which we infer that $h(G(a))\subseteq G^{c}(h(a))$. Conversely, assume that 
 $F\in G^{c}(h(a))$. Then for every $P\in{\mathcal X}(W)$,  $FR^{c}P$ implies that $P\in h(a)$. Suppose that $G(a)\notin F$. Then $G^{-1}(F)$ is a filter and $a\notin G^{-1}(F)$. Hence, there is $T\in {\mathcal X}(W)$ such that $a\notin T$ and $G^{-1}(F)\subseteq T$. This last assertion and (F2) allow us to conclude that $FR^{c}T$. From this statement we have that $T\in h(a)$ and so, $a\in T$, which is a contradiction. Therefore, $h(G(a))=G^c(h(a))$. Similarly, it is shown that $h(H(a))=H^c(h(a))$. Thus, by virtue of the results established in \cite{EO3, EO} the proof is completed.

Lemma \ref{l23} will show that the order-embedding $k:X\to {\mathcal X}({\mathcal C}(X))$ defined by  $k(x)=\{A\in{\mathcal C}(X): x\in A\}$ for every $x\in X$ (see \cite{EO3, EO} ) preserves the relations $R$ and $Q$.

\begin{lemma}\label{l23}
Let  $(X,\leq,g,R,Q)$ be a $TSH$-frame and let $x,y\in X$. Then
\begin{itemize}
\item [\rm (i)] $xRy$ if and only if $k(x)R^{c}k(y)$,
\item [\rm (ii)] $xQy$ if and only if $k(x)Q^{c}k(y)$.
\end{itemize}
\end{lemma}

{\bf Proof}

We will only prove (i). Assume that $xRy$ and suppose that $A\in{\mathcal C}(X)$ verifies $G^c(A)\in k(x)$. Then it is easy to see that  $y\in A$ and so, $k(x)R^{c}k(y)$. Conversely, let $x,y\in X$ be such that $k(x)R^{c}k(y)$. Then ${G^c}^{-1}(k(x))\subseteq k(y)$. On the other hand, note that $[\leq](X\setminus (y])\in{\mathcal C}(X)$ and $y\notin [\leq ](X\setminus (y])$. Thus, $[\leq](X\setminus (y])\notin k(y)$ and so, $[\leq](X\setminus (y])\notin {G^c}^{-1}(k(x))$. Therefore, $[R]([\leq](X\setminus (y]))\notin k(x)$ from which we infer that $x\notin[R]([\leq](X\setminus (y]))$. Then there is $z$ such that $xRz$ and $z\notin [\leq](X\setminus(y])$. From this last assertion there is $w$ such that $z\leq w$ and $w\leq y$, which allow us to infer that $z\leq y$. Hence, by virtue of the reflexivity of $\leq$ and (K3), $xRy$ as required.

Lemma \ref{l23} and the results indicated in \cite{EO3, EO} enable us to conclude

\begin{theorem}\label{t22}
Every $TSH$-frame $X$ is embeddable into the canonical frame of its complex algebra ${\mathcal X}({\mathcal C}(X))$. 
\end{theorem}

Theorems \ref{t21} and \ref{t22} enable us to obtain a discrete duality for $TSH$-algebras.

\section{\large A propositional calculus based on $TSH$-algebras}

In this section, we will describe a propositional calculus that has $TSH$-algebras as the algebraic counterpart. The terminology and symbols used here coincide in general with those used in \cite{Ra}.

\vspace{1mm}

Let ${\mathcal L}=(A^0,For[V])$ be a formalized language of zero order, where in the alphabet $A^0=(V, L_0,L_1,L_2,U)$ the set

\begin{itemize}
\item $V$ of propositional variables is enumerable,
\item $L_0$ is empty,
\item $L_1$ contains three elements denoted by  $\sim$, $G$ and $H$ called negation sign and tense operators signs, respectively,
\item $L_2$ contains three elements denoted by $\vee$,\, $\wedge$,\, $\to$,\, called disjunction sign, conjunction sign and implication sign, respectively,
\item $U$ contains two elements denoted by $($, $)$.
\end{itemize}

\vspace{1mm}

For any $\alpha,\, \beta$ in the set $For[V]$ of all formulas over $A^0$, instead of $(\alpha\to \beta)\wedge (\beta\to \alpha)$, $\sim G\sim \alpha$ and $\sim H\sim \alpha$ we will write for brevity $\alpha\leftrightarrow \beta$, $F\alpha$ and $P\alpha$, respectively.

\vspace{1mm}

We assume that the set ${\mathcal A}_l$ of logical axioms consists of all formulas of the following form, where $\alpha,\, \beta, \gamma$ are any formulas in $For[V]$:

\begin{itemize}
\item [(M0)] the axioms of the symmetric modal propositional calculus, i.e., the axioms (A1)-(A10) indicated in  \cite[page 60]{M},
\item[(M1)] $G(\alpha \to \beta) \to (G\alpha\to G\beta )$,\,   $H(\alpha\to \beta )\to (H\alpha\to  H\beta)$,
\item[(M2)] $\alpha \to GP\alpha$,\,  $\alpha\to HF\alpha$.
\end{itemize}

\vspace{1mm}

The consequence operation $C_{\mathcal L}$ in ${\mathcal L}$ is determined by ${\mathcal A}_l$ and by the following rules of inference:

\vspace{1mm}

\begin{tabular}{ll}
\hspace{-.1cm}(R1) $\displaystyle\frac{\alpha,\;\; \alpha\rightarrow \beta}{\beta}$, &  \hspace{2.4cm} (R3)  $\displaystyle\frac{\alpha}{G\alpha}$,\\
\hspace{-.1cm}(R2)  $\displaystyle\frac{\alpha\to\beta}{\sim\beta\rightarrow\, \sim \alpha}$, & \hspace{2.4cm} (R4)  $\displaystyle\frac{\alpha}{H\alpha}$.
\end{tabular}

\vspace{1mm}

The system $\mathcal{TMS}=({\mathcal L}, C_{\mathcal L})$ thus obtained will be called the $\mathcal {TMS}$--proposi\-tional calculus. We will denote by ${\mathcal T}$ the set of all formulas derivable in $\mathcal {TMS}$. If $\alpha$ belongs to ${\mathcal T}$ we will write  $\vdash\alpha$.

\vspace{1mm}

Let $\approx$ be  the binary relation on $For[V]$ defined by

\begin{center}
$\alpha\approx \beta$ if and only if $\vdash\alpha\leftrightarrow \beta$.
\end{center}

Then it is easy to check that  $\approx$ is a congruence relation on $\langle For[V], \vee,\wedge,\to,\sim,G,H\rangle$ and ${\mathcal T}$ determines an equivalence class which we will denote by $1$. Moreover, taking into account \cite[page 62]{M} it is straightforward to prove

\begin{theorem} $\langle For[V]/\approx,\vee,\wedge,\to,\sim, G, H, 0,1\rangle$ is a $TSH$-algebra, being $0=\sim 1$.
\end{theorem}

\begin{definition} A $TSH$-model based on a $TSH$-frame $K=(X,\leq,g,R,Q)$ is a system $M=(K,m)$ such that $m:V\to {\mathcal P}(X)$ is a meaning function that assigns subsets of states to propositional variables, i.e. satisfies the following condition:
\begin{itemize}
\item[\rm (her)] $x\leq y$ and $x\in m(p)$ imply $y\in m(p)$.
\end{itemize}
\end{definition}

\vspace{1mm}

\begin{definition} 
A $TSH$-model $M=((X,\leq,g,R,Q); m)$  satisfies a formula $\alpha$ at the  state $x$ and we write $M\models_{x}\alpha$, if  the following conditions are satisfied:

\begin{itemize}
\item $M\models_{x}p$ if and only if $x\in m(p)$ for $p\in V$,
\item $M\models_{x}\alpha\vee \beta$ if and only if $M\models_{x}\alpha$ or $M\models_{x}\beta$,
\item $M\models_{x}\alpha\wedge \beta$ if and only if $M\models_{x}\alpha$ and $M\models_{x}\beta$,
\item $M\models_{x}\sim\alpha$ if and only if $M\not\models_{g(x)}\alpha$,
\item $M\models_{x}\alpha\to \beta$ if and only if for all $y$, if $x\leq y$ and $M\models_{y}\alpha$ then $M\models_{y}\beta$,
\item $M\models_{x}G\alpha$ if and only if for all $y$, if $xRy$ then $M\models_{y}\alpha$,
\item $M\models_{x}H\alpha$ if and only if for all $y$, if $xQy$ then $M\models_{y}\alpha$.
\end{itemize}
\end{definition}

A formula $\alpha$ is {\em true in a $TSH$-model} $M$ (denoted by $M\models\alpha$) if and only if for every $x\in W$, $M\models_x\alpha$.
The formula $\alpha$ is {\em true in a $TSH$-frame} $K$ (denoted by $K\models\alpha$) if and only if it is true in every $TSH$-model based on $K$. The formula $\alpha$ is $TSH$-{\em valid} if and only if it is true in every $TSH$-frame.

\begin{proposition} Given a $TSH$-model $M=((X,\leq,g,R,Q); m)$, the meaning function $m$ can be extended to all formulae by $m(\alpha)=\{x\in X: M\models_x\alpha\}$. For every $TSH$-model $M$ and for every formula $\alpha$, this extension has the property
\begin{itemize}
\item [\rm (her)] if $x\leq y$ and $x\in m(\alpha)$ then $y\in m(\alpha)$.
\end{itemize}
\end{proposition}

{\bf Proof} The proof is by induction with respect to complexity of $\alpha$. By way of an example we show ({\rm her}) for formulas of the form $G\alpha$. Let (1) $x\leq y$ and (2) $M\models_{x}G(\alpha)$. Suppose that $yRz$, then by (1),(2) and (K3), we have $M\models_{z}\alpha$.

\begin{theorem}\label{t32}{\bf (Completeness Theorem)} Let $\alpha$ be a formula in $\mathcal{TMS}$. Then the following conditions are equivalent:
\begin{itemize}
\item[\rm (i)] $\alpha$ is derivable in $\mathcal {TMS}$,
\item[\rm (ii)] $\alpha$ is $TSH$-valid. 
\end{itemize}\end{theorem}

{\bf Proof} ${\rm (i)}\Rightarrow{\rm(ii)}$: We proceed by induction on the complexity of the formula $\alpha$. For
example, we shall prove that the axiom (M2) is $TSH$-valid.
Let $K=(X,\leq,g,R,Q)$ be a $TSH$-frame and $M$ a $TSH$--model based on $K$.
\begin{itemize}
\item[(1)] Let $y\in X$ be such that $x\leq y$,\hfill[hip.]
\item[(2)] ${M}\models_{x}\alpha$,\hfill[hip.] 
\item[(3)] Let $z\in X$ be such that $y\, Q\, z$,\hfill[hip.]  
\end{itemize}
Suppose that
\begin{itemize}
\item[(4)] ${M}\models_{g(z)}G\sim \alpha$,\hfill[hip.]
\item[(5)] $x\, Q\, z$,\hfill[(1),(3), (K3)]  
\item[(6)] $g(z)\, R\, g(x)$,\hfill[(5),(K5)]
\item[(7)] ${M}\models_{g(x)}\sim \alpha$,\hfill[(4),(6)] 
\item[(8)] ${M}\not\models_{x}\alpha$. \hfill[(7),(K2)] 
\end{itemize}

(8) contradicts (2). Then

\begin{itemize}
\item[(9)] ${M}\not\models_{g(z)}G\sim \alpha$,\hfill[(4),(8)]
\item[(10)] ${M}\models_{z}\sim G\sim \alpha$,\hfill[(9)]
\item[(11)] ${M}\models_{z}H\sim G\sim \alpha$,\hfill[(3),(10)]
\item[(12)]${M}\models_{x}\alpha\to H\sim G\sim \alpha$.\hfill[(1),(2),(11)]
\end{itemize}

\noindent ${\rm (ii)}\Rightarrow{\rm(i)}$:  Assume that $\alpha$ is not derivable, i.e. $[\alpha]_\approx \not= 1$. We apply Theorem \ref{t21} to
the $TSH$-algebra $For[V]/\approx$, hence there exists a $TSH$-frame ${\mathcal X}(For[V]/\approx)$ and an injective
morphism of $TSH$-algebras $h:For[V]/\approx\rightarrow {\mathcal C}({\mathcal X}(For[V]/\approx))$. Let us consider the function $m:\mathcal {TMS}\to{\mathcal C}({\mathcal X}(For[V]/\approx))$ defined by $m(\alpha)=h([\alpha]_\approx)$ for all $\alpha \in For[V]$. It is straightforward to prove that $m$ is an meaning function. Since $h$ is injective, $m(\alpha)=h([\alpha]_\approx)\not={\mathcal X}(For[V]/\approx)$, i.e. $({\mathcal X}(For[V]/\approx),m)\not\models_{x_o}\alpha$ for some $x_o\in {\mathcal X}(For[V]/\approx)$.  Thus $\alpha$ is not $TSH$-valid.


\begin{thebibliography}{99}

\bibitem[1]{EO3} W. Dzik, E. Orlowska and C. van Alten,
           {\it Relational representation theorems for general lattices with negations},
            Relations and Kleene algebra in computer science, 162--176, Lecture Notes in Comput. Sci., 4136, Springer, Berlin, 2006.

\bibitem[2]{FPS} A.V. Figallo, G. Pelaitay, C. Sanza,
             {\it Operadores temporales sobre \'algebras de Heyting sim\'etricas},
              Noticiero de la Uni\'on Matem\'atica Argentina (2009).


\bibitem[3]{I2}L. Iturrioz, {\it Symmetrical Heyting algebras with operators},
           Z. Math. Logik Grundlag. Math. {\bf 29} (1983), no.1, 33--70. 

\bibitem[4]{K} T. Kowalski,
           {\it Varieties of tense algebras},
           Rep. Math. Logic. {\bf 32} (1998), 53--95.  

\bibitem[5]{Mo} Gr. C. Moisil,
           {\it Logique modale},
           Disquisit. Math. Phys. {\bf 2} (1942), 3--98. 


\bibitem[6]{M} A. Monteiro,
          {\it Sur les alg\`ebres de Heyting Simetriques},
           Special issue in honor of Ant\'onio Monteiro. Portugal. Math. {\bf 39}(1--4), 1980.


\bibitem[7]{EO} E. Or\l owska and I. Rewitzky,
           {\it Discrete Dualities for Heyting Algebras with Operators}, 
            Fund. Inform. {\bf 81} (2007), no.1--3, 275--295.

\bibitem[8]{EO1} E. Or\l owska and I. Rewitzky, 
          {\it Duality via Truth: Semantic frameworks for lattice-based logics},
           Log. J. IGPL {\bf 13} (2005), no. 4, 467--490.


\bibitem[9]{Ra} H. Rasiowa,
           {\it An algebraic approach to non-classical logics},
            Studies in Logic and the Foundations of Mathematics, 78, North-Holland, 1974.

\bibitem[10]{S} H.P. Sankappanavar,
           {\it Heyting algebras with a dual lattice endomorphism},
            Z. Math. Logik Grundlag. Math. {\bf 33} (1987), no. 6, 565--573.

\bibitem[11]{PU} P. Unterholzner,
            {\it Algebraic and relational semantics for tense logics}.
             Rend. Sem. Mat. Univ. Padova {\bf 65} (1981), 119--128.


\end{thebibliography}
\end{document}